\documentclass{article}

\newtheorem{Th}{Theorem}

\newtheorem{Lemma}{Lemma}
\newtheorem{Cor}{Corollary}
\newtheorem{Def}{Definition}
\newtheorem{Rem}{Remark}
\newcommand{\weg}[1]{} 
\newcommand{\wegg}[1]{}
\newcommand{\eqdef}{\stackrel{{\rm def}}{=}}
\newcommand{\De}{\Delta}
\newcommand{\de}{\delta}

\newcommand{\si}{\sigma}
\newcommand{\om}{\omega}
\newcommand{\bg}{\bar g}

\newcommand{\dif}{\mbox{\rm d}}

\title{Geodesic equivalence and integrability.}
\author{  Peter J.  Topalov, Vladimir S. Matveev}
\date{}
\begin{document}
\maketitle


{
 \begin{abstract}
  We suggest a construction that, given a 
 trajectorial diffeomorphism between   
 two Hamiltonian systems,  produces integrals of them.  
 As the main  example we treat geodesic equivalence of metrics.  
 We show that the  existence of a non-trivially geodesically 
 equivalent metric leads to Liouville integrability, and present  
 explicit formulae for integrals. 
\end{abstract}
}

\section{Introduction} 
Integrals of a system are closely  related  to symmetries. 
A classical example is   Noether's theorem:
if  a vector field $X$ on a manifold  $M$  preserves a 
Lagrangian $L:{\cal   T}M\to R$ 
, then the function
$I_X\eqdef\frac{\partial L}{\partial{\dot x}}(x,\dot x)X(x)$ is a first
integral of the corresponding Lagrangian system. 

There are many generalizations of Noether's  theorem, 
we recall the following two. 
In the paper \cite{Koz-Bol} it was shown  that the  existence of 
a vector field  on ${\cal T}^*M$ which  commutes  with  
a Hamiltonian vector field 
 allows one  to construct a (multi-valued) integral of the 
Hamiltonian system. In the paper \cite{TensInv} the result of
\cite{Koz-Bol}  was generalized to 
tensor fields. It was shown that if a Hamiltonian flow 
 preserves  a tensor field on   ${\cal T}^*M$, then  there  exists  an  
(also multi-valued) integral of the Hamiltonian system.

In our   paper we,   following ideas of \cite{TensInv},
 present a construction  which, 
given  
   a  diffeomorphism between two Hamiltonian systems 
   that 
   takes the  trajectories  and the isoenergy surfaces
  of the first Hamiltonian system    
 to the trajectories and  the  isoenergy surfaces of 
the second one, produces  $n$ integrals  
of the first system, where $n$ is the number of 
the degrees of freedom of the system.

The construction is applied to geodesically equivalent metrics. 
Let $g=(g_{ij})$ and  $\bar g=(\bar g_{ij})$ be  smooth 
metrics on the  same manifold $M^n$. 
\begin{Def}\hspace{-2mm}{\bf .}  The 
metrics $g$  and $\bar g$  
are {\em geodesically equivalent},
if they have  the same geodesics  (considered as unparameterized curves).
\end{Def}

 This is  rather classical material. In  1869 Dini \cite{Dini} formulated the problem of 
local  classification of geodesically equivalent metrics, and solved it for 
 dimension two.  
 In 1896 Levi-Civita \cite{Levi-Civita}
 got a  local description of geodesically
equivalent metrics on manifolds of arbitrary dimension.
In the paper \cite{Mikes} a family of (non-trivial) examples of geodesically 
equivalent metrics on closed manifolds was constructed.

For geodesically equivalent metrics, a  
trajectorial diffeomorphism 
 $\Phi$ is given by 
$\Phi(x,\xi)=(x, \frac{\|\xi\|_g}{\|\xi\|_{\bar g}}\xi)$.
Here $(x,\xi)\in {\cal T}M^n$, $x$ is a point  of $M^n$ and $\xi\in {\cal T}_xM^n$.

\begin{Th}\hspace{-2mm}{\bf .} \label{Main_Topalov}
Let metrics  $g$ and $\bar g$ on $M^n$  
be geodesically equivalent. 
Denote by $G$ the linear operator $g^{-1}\bar g=(g^{i\alpha}\bar g_{\alpha j})$. Consider the characteristic polynomial  
$det(G- \mu E)=c_0\mu^n+c_1\mu^{n-1} + ... + c_n$. The
coefficients $c_1, .. , c_n$ are smooth functions on the 
manifold $M^n$, and $c_0\equiv (-1)^n$. 
Then 
the functions 
$I_k=\left(\frac{det(g)}{det(\bar g)}\right)^{\frac{k+2}{n+1}}\bar g(S_k
\xi,\xi)$, $k=0, ... , n-1$, where $S_k\eqdef\sum_{i=0}^kc_i G^{k-i}
$, are
 integrals of the geodesic flow of the metric $g$
and pairwise  commute.  
\end{Th}

\begin{Rem}\hspace{-2mm}{\bf .}   The integral $I_0=\left(\frac{det(g)}{det(\bar g)}\right)^{\frac{2}{n+1}}\bar g(\xi, \xi)$  was obtained by   
Painlev\'e, see \cite{Levi-Civita}. The integral 
 $I_{n-1}$  is  the energy integral (multiplied by minus two).  
\end{Rem}
The integrals $I_1, I_{2},...,I_{n-2}$  seem to be new, although 
 in each Levi-Civita chart the integrals are linear combinations 
of  Levi-Civita integrals
(see Section~\ref{LCtheorem} for definitions). 
We touch on 
the connection between the integrals $I_0,..., I_{n-1}$
  and Levi-Civita integrals  
in Section \ref{LiouvIntegr}.

 Metrics $g, \bar g$ on $M^n$ 
are {\em strictly non-proportional} at  a point $x\in M^n$, if 
the characteristic polynomial 
$\frac{1}{det(g)}det(\bar g-tg))_{|x}$ 
has no multiple  root. 

\begin{Cor}\hspace{-2mm}{\bf .} \label{T} Let $M^n$ be a closed real-analytic 
manifold supplied with
  two real-analytic metrics $g,\bar g$ such that the metrics $g,\bar g$ are 
 geodesically equivalent and 
                 strictly non-proportional at  least at one point.     
Then  the fundamental group 
$\pi_1(M^n)$ of the manifold $M^n$ contains  a commutative   subgroup of finite 
index, and the dimension of the homology  group $H_1(M^n; {\bf Q})$  is 
 no greater than $n$.  
\end{Cor}

For dimension two the  converse of Theorem~\ref{Main_Topalov} is also true, 
and  the condition   
of Corollary~\ref{T} can be weakened.

\begin{Cor}\hspace{-2mm}{\bf .} \label{dim2}
Metrics  $g$ and $\bar g$ on a surface $M^2$  are geodesically equivalent, 
if and only if the function $\left(\frac{det(g)}{det(\bar g)}\right)^{\frac{2}{3}}\bar g(\xi, \xi) $ is an integral of the geodesic flow of the metric 
$g$. 
\end{Cor}

\begin{Cor}\hspace{-2mm}{\bf .}{\label{Kol}} Let metrics $g, \bar g$ on a closed 
surface of 
negative Euler characteristic be 
geodesically equivalent.
 Then $g=C  \bar g$, where $C$ is a constant.  
\end{Cor}

\begin{Cor}\hspace{-2mm}{\bf .}
\label{torus}
Let metrics $g, \bar g$  
  on the torus $T^2$
 be geodesically equivalent. 
If they  are proportional at  
a point $x\in T^2$, then $g=C  \bar g$, where $C$ is a positive constant.  
\end{Cor}

\begin{Cor}\hspace{-2mm}{\bf .}\label{sphere}
Let metrics $g$, 
$\bar g$ on the sphere $S^2$  be geodesically equivalent. 
Then there are three possibilities. 
\begin{itemize}
\item[1.] The metrics are proportional  at   exactly   two points.
\item[2.] The metrics are proportional  at  exactly   four points. 
\item[3.] The metrics are completely proportional, i.e. $g=C  \bar g$,
 where $C$ is a positive constant.   
\end{itemize}
In the first case the  metrics admit a Killing vector field. 
\end{Cor}

Recall that a vector field on $M^n$  is 
{\it Killing} (with respect to a metric), 
 if the flow of the field 
 preserves the metric.

\begin{Cor}\hspace{-2mm}{\bf .}\label{nontrivial}
Let  metrics $g, \bar g$  on a surface $M^2$ be geodesically equivalent.
If the metrics are proportional at each point of an open non-empty 
domain $U\subset M^2$, then  $g=C  \bar g$, where $C$ is a positive constant.   \end{Cor}

\begin{Cor}\hspace{-2mm}{\bf .}\label{killing}
If metrics $g,\bar g$ on a manifold  $M^n$ are geodesically equivalent, and if
 the 
metric $g$ admits  a non-trivial  Killing vector field, then the metric $\bar g$ also admits a non-trivial Killing vector 
field. 
\end{Cor}

One of the most famous integrable geodesic flows on closed surfaces
is the geodesic flow of the metric on ellipsoid (see \cite{Moser}). 
Consider the ellipsoid \newline $\sum_{i=1}^n\frac{(x^i)^2}{a_i}=1$, where $a_i>0$, $i=1,...,n$.  
\begin{Th}\hspace{-2mm}{\bf .}\label{ellipsoid}
The restriction of the metric $\sum_{i=1}^n(dx^i)^2$ to 
the ellipsoid 
 \newline $\sum_{i=1}^n\frac{(x^i)^2}{a_i}=1$ is geodesically equivalent to 
the restriction of the metric 
$$
 \frac{1}{\sum_{i=1}^{n}\left(\frac{x^i}{a_i}\right)^2}\left(\sum_{i=1}^n\frac{(dx^i)^2}{a_i}\right)
$$
to the ellipsoid. 
\end{Th}

The paper is organized as follows. In Section \ref{Theory} we present the 
announced construction. Theorem \ref{MainTh} there gives an explicit  
formula for  a
one-parameter family of first integrals,  if a trajectorial
diffeomorphism between two Hamiltonian systems is given.

In Section \ref{LCtheorem}, for use in Sections \ref{GeodEqInt},
\ref{LiouvIntegr}, \ref{Ellipsoid} 
we formulate 
 Levi-Civita and Painlev\'e  results about a local form 
of geodesically equivalent  metrics.

In Section  \ref{GeodEqInt}  
we apply the construction 
 to   geodesically equivalent
metrics, and prove   that the functions  $I_0, ..., I_{n-1}$  from 
  Theorem \ref{Main_Topalov} are integrals of the geodesic flow of the metric $g$. 

In Section \ref{LiouvIntegr} we prove that the integrals  $I_0, ..., I_{n-1}$ 
are in involution.

In Section \ref{Obstructions} we prove  
Corollaries \ref{T}, \ref{dim2},  \ref{Kol}, \ref{torus}, \ref{sphere}, \ref{nontrivial},  \ref{killing}.

In Section \ref{Ellipsoid} we prove Theorem~\ref{ellipsoid}.

The authors are grateful to 
A.~V.~Bolsinov, A.~T.~Fomenko, 
V.~V.~Kozlov  and 
I.~A.~Taimanov  for useful discussions. The main   results of the paper
were
obtained during a 4-week visit of P.~Topalov to Bremen University.
 The authors are
grateful to the Institute of Theoretical  Physics of  Bremen University
 for its  hospitality and to
 the Deutsche Forschungsgemeinschaft for  partial financial support.

\section{Trajectorial diffeomorphisms  and integrals}\label{Theory}

\noindent Let  $v$ and $\bar v$ be Hamiltonian systems
 on  
 symplectic manifolds $(M,\omega)$ and $(\bar M,\bar\omega)$
 with Hamiltonians $H$ and
$\bar H$ respectively. 
Consider the 
isoenergy surfaces  
$$
Q\eqdef\left\{x\in M: H(x)=h\right\}, \ \ 
\bar Q\eqdef\left\{x\in \bar  M: \bar H(x)=\bar h\right\},
$$
\noindent where $h$ and $\bar h$ are regular values of the functions 
$H$, $\bar H$ respectively. 
Let $U(Q)\subset M$ and  $U(\bar Q)\subset \bar M$ 
be neighborhoods of the isoenergy surfaces $Q$ and $\bar Q$.

\begin{Def}\hspace{-2mm}{\bf .}\label{TrAutDef}
A diffeomorphism $\Phi : U(Q)\longrightarrow U(\bar Q)$, $\Phi(Q)=\bar Q$, 
 is said to be  {\em
trajectorial on $Q$}, 
 if the restriction $\Phi|_Q$
 takes the trajectories of the system $v$ to
the trajectories of the system $\bar v$.
\end{Def}

Denote the restriction $\Phi|_Q$ by $\phi$. Since $\phi$ takes 
the trajectories of $v$ to the trajectories of $\bar v$, it takes 
the vector field $v$ to the vector field that is proportional to 
$\bar v$. Denote by $a_1:Q\to R$ 
the coefficient of proportionality, i.e. $\phi_*(v)=a_1\bar v$. 
Since $\Phi$ takes 
 $Q$ to  $\bar Q$, 
it takes 
the differential  $\dif H$ 
to a form  that is proportional
 to 
$d\bar H$. Denote by $a_2:Q\to R$ 
the coefficient of proportionality, i.e. $\phi_* \dif H=a_2\dif \bar H$. 
By $a$ we denote the product  $a_1a_2$. 
We denote the Pfaffian of a skew-symmetric matrix $X$ by Pf($X$).

\begin{Th}\hspace{-2mm}{\bf .}\label{MainTh}
Let a diffeomorphism 
$\Phi:{U}(Q)\to{U}(\bar Q)$, $\Phi(Q)=\bar Q$, 
 be trajectorial on $Q$.
 Then for each value of the parameter $t$ 
 the polynomial
\begin{equation}\nonumber
{\cal P}^{n-1}(t)\eqdef
\frac{\mathop{\rm Pf}\left(\Phi^{*}\bar\om-t\om\right)}
{\mathop{\rm Pf}\left(\om\right)(t-a)}
\end{equation}
is an integral of the system $v$ on $Q$.
 In particular, all the
coefficients of the polynomial ${\cal P}^{n-1}(t)$  are 
 integrals.
\end{Th}

\noindent {\bf Proof.} 
Denote by $\sigma$, $ \bar \sigma$
the restrictions of the forms $\omega, \bar \omega$  to $Q$, $\bar Q$  
respectively.    
Consider the form $\phi^*\bar \sigma$ on $Q$.

\begin{Lemma}[Topalov, \cite{TensInv}]\hspace{-2mm}{\bf .}\label{TensInv}
The flow  $v$ preserves the
 form  $\phi^{*}\bar\sigma$.
\end{Lemma}

\noindent {\bf Proof of Lemma \ref{TensInv}.}
The  Lie derivative $L_v$ of the form $\phi^*\bar \sigma$ along the vector field $v$ satisfies  
$$
L_v\phi^{*}\bar\sigma=\dif\left[{\imath}_v\phi^{*}\bar\sigma\right]+{\imath}_v
\dif\left[\phi^{*}\bar\sigma\right]. 
$$
On the right side  both terms vanish. 
More precisely, for an arbitrary  vector
$u\in {\cal T}_xQ$ at an arbitrary point $x\in Q$ 
we have    
$$
\begin{array}{rcl}
{\imath}_{v}\phi^{*}\bar\sigma(u)&= \ \bar\sigma(\phi_{*}(v),\phi_{*}(u)) & =\\
                     &= \ \bar\sigma(a_1\bar v,\phi_{*}(u)) & =\\
                     &= \  -a_1\dif\bar{\it H}(\phi_{*}(u)) & =0. 
\end{array}
$$
\noindent Since the form $\bar\omega$ is closed, the form $\bar \sigma$ is also closed and
$\dif\left[\phi^{*}\bar\sigma\right]=\phi^{*}(\dif\bar\sigma)=0$, q.~e.~d.

It is obvious  that the kernels of the forms $\sigma$ and
$\phi^{*}\bar\sigma$ coincide (in the space ${\cal T}_xQ$ at each point $x\in Q$) 
with the linear  span of
the  vector  $v$. Therefore these forms induce two
non-degenerate
tensor fields on the quotient bundle 
${\cal T}Q/\langle v\rangle$. We shall
denote the corresponding forms on ${\cal T}Q/\langle v\rangle$ also by the  letters $\sigma, \bar\sigma$.

\begin{Lemma}\hspace{-2mm}{\bf .}\label{Din}
The characteristic polynomial of the operator
$(\sigma)^{-1}(\phi^{*}\bar\sigma)$ on ${\cal T}Q/\langle v\rangle$
is preserved by the flow $v$.
\end{Lemma}

\noindent{\bf Proof of Lemma~\ref{Din}.} Since 
 the flow $v$ preserves the 
Hamiltonian  $H$ and the form $\omega$,   
the flow $v$ preserves 
the form $\sigma$. 
Since the  flow  $v$ preserves 
both forms, it preserves  
the  characteristic
polynomial of the  operator 
$(\sigma)^{-1}({\phi^{*}\bar\sigma})$, q.~e.~d.

Since both forms are skew-symmetric, each root of the  characteristic
 polynomial 
of  the operator
$(\sigma)^{-1}(\Phi^{*}\bar\sigma)$ has an even multiplicity. Then  the characteristic
 polynomial   is the  square of a polynomial 
$\delta^{n-1}(t)$ of degree $n-1$. Hence  
the polynomial $\delta^{n-1}(t)$ is also preserved by the flow $v$.
It is obvious  that
\begin{equation}\nonumber
\delta^{n-1}(t)=(-1)^{n-1}
\frac{\mathop{\rm Pf}\left(\phi^{*}\bar\sigma-t\sigma\right)}
{\mathop{\rm Pf}\left(\sigma\right)}.
\end{equation}

The last step of the  proof is to verify
that 
$$(t-a)\delta^{n-1}=\frac{\mathop{\rm Pf}\left(\Phi^{*}\bar\omega-t\omega\right)}
{\mathop{\rm Pf}\left(\omega\right)}\eqdef\Delta^n.$$ Take an arbitrary point $x\in Q$. 
Consider the  form $\Phi^*\bar\omega-a\omega$ on 
${\cal T}_xM$.
The form $\imath_v(\Phi^*\bar\omega-a\omega)$ equals zero. 
More precisely, 
 for  any  vector  $u\in{\cal T}_xM$  we have 
\begin{eqnarray*}
\imath_{v}(\Phi^{*}\bar\om-a\om)&=&\bar\om(\Phi_*(v),\Phi_*(u))-a\om(v,u)=\\
                                 &=&\bar\omega(a_1v, \Phi_*(u))-a\om(v,u)=\\
                                 &=&-a_1\dif\bar H(\Phi_*(u))+adH=\\
                                 &=&-a\dif H+adH=0.
\end{eqnarray*}
\noindent  There exists a vector $A\in {\cal T}_xM$ such that 
$\om(A,v)\ne 0$ and the restriction 
of the form $\imath_A(\Phi^*\bar\omega-a\omega)$ to the space   ${\cal T}_xM$
 equals zero. 
More precisely, since the forms 
 $\Phi^*\bar\omega$, $\omega$ are skew-symmetric, then the kernel 
$K_{\Phi^*\bar\omega-a\omega}$ of the form 
$\Phi^*\bar\omega-a\omega$ has an even dimension, and the  
kernel 
of the restriction of the form $\Phi^*\bar\omega-a\omega$ to 
${\cal T}_xQ$ has an odd dimension. 
Thus  the intersection  
 $ K_{\Phi^*\bar\omega-a\omega}\cap ({\cal T}_xM\setminus {\cal T}_xQ)$ 
is not empty. For 
 each vector $A$ from the intersection  we obviously have  
$\om(A,v)\ne 0$ and $\imath_A(\Phi^*\bar\omega-a\omega)=0$. 
Without loss of generality 
 we can   assume   $\om(A,v)=1$.

 Consider   a basis $(v,e_1,...,e_{2n-2})$ for  the space 
${\cal T}_xQ$. 
The set  $(A,v,e_1,...,e_{2n-2})$ is a basis 
 for  the space 
${\cal T}_xM$. In this basis we have  
$$
\begin{array}{rcl}
\det(\Phi^{*}\bar\om-t\om)&=&
\det\left|
\begin{array}{cc|c}
0&a-t&(*)\\
-(a-t)&0&0\cdots 0\\
\hline

-(*)&0&(\Phi^{*}\bar\om-t\om)_{\langle e_1,...,e_{2n-2}\rangle}\\
\end{array}
\right|\\ & & \\
&=&(a-t)^2\det((\Phi^{*}\bar\om-t\om)_{\langle e_1,...,e_{2n-2}\rangle})\\
&=&(a-t)^2\det(\phi^*\bar\sigma-t\sigma),
\end{array}
$$
\noindent where 
$(\Phi^{*}\bar\om-t\om)_{\langle e_1,...,e_{2n-2}\rangle}$  is
 the
matrix of the form $\Phi^{*}\bar\om-t\om$ in the basis 
$(e_1,...e_{2n-2})$. Finally,   $\delta^{n-1}= {\cal P}^{n-1
}$,  q.~e.~d .

\section{Levi-Civita theorem}\label{LCtheorem}

\noindent Let $g$ and $\bar g$ be smooth metrics  on a manifold $M^{n}$.
Recall that the  common eigenvalues 
of the metrics $g,\ \bar g$ at a
point $x\in M$ are roots of the
characteristic polynomial $P_x(t)={\rm  det}\left(G-tE\right)_{|x}$, where 
$G\eqdef\left(g^{i\alpha }\bar g_{\alpha j}\right)$.
Suppose that
 at every   point of
 an open domain ${\cal D}\subset M^n$
 the common eigenvalues of the metrics $g, \ \bar g$ assume
 $m$ distinct  values
$\rho^1,\rho^2,...,\rho^m$ \  $(1\le m\le n)$  with  multiplicities
 $k_1,k_2,...,k_m,$
respectively.

In the  paper \cite{Levi-Civita},
 Levi-Civita 
 proved 
that for every point
$P\in{\cal D}$ there is 
 an open neighborhood ${\cal U}(P)\subset{\cal D}$ and
a coordinate system
$\bar x=(\bar x_1,...,\bar x_m)$
(in ${\cal U}(P)$), where $\bar x_i=(x_i^1,...,x_i^{k_i})$,
$(1\le i\le m)$, such that the quadratic  forms  of the 
metrics $g$ and $\bg$ have the following form: 
\begin{eqnarray}
g(\dot{\bar x}, \dot{\bar x})&=&\ \Pi_1(\bar x)A_1(\bar x_1,\dot{\bar x}_1)+
\ \Pi_2(\bar x)A_2(\bar x_2,\dot{\bar x}_2)+\cdots+\nonumber\\
&+&\ \Pi_m(\bar x)A_m(\bar x_m,\dot{\bar x}_m),\label{Canon_g}\\
\bg(\dot{\bar x}, \dot{\bar x})&=&\rho^1\Pi_1(\bar x)A_1(\bar x_1,\dot{\bar x}_1)+
\rho^2\Pi_2(\bar x)A_2(\bar x_2,\dot{\bar x}_2)+\cdots+\nonumber\\
&+&\rho^m\Pi_m(\bar x)A_m(\bar x_m,\dot{\bar x}_m)\label{Canon_bg},
\end{eqnarray}

\noindent where $A_i(\bar x_i, \dot{\bar x}_i)$ are  
 positive-definite quadratic forms  in  the velocities $\dot{\bar x}_i$
with coefficients depending  on $\bar x_i$,  
\begin{eqnarray}
\Pi_i&\eqdef&(\phi_i-\phi_1)\cdots(\phi_{i}-\phi_{i-1})(\phi_{i+1}-\phi_i)\cdots
(\phi_m-\phi_i)
\end{eqnarray}
and  $\phi_1,\phi_2,..., \phi_m$, $0<\phi_1<\phi_2<...<\phi_m$,  are smooth functions such that
$$
\phi_i=\left\{
\begin{array}{l}
\phi_i(\bar x_i),\mbox{\quad if}\quad k_i=1\\
{\rm constant},\quad\mbox{else}.
\end{array}
\right.
$$
It is easy to see that the functions $\rho^i$ as functions of $\phi_i$
and the function $\phi_i$ as  functions of $\rho^i$   are
given by 
\begin{eqnarray}
  \rho^i  & = & \frac{1}{\phi_1...\phi_m}\frac{1}{\phi_i} \nonumber \\
  \phi_i  & = & \frac{1}{\rho^i}(\rho^1\rho^2...\rho^m)^{\frac{1}{m+1}}
   \nonumber
\end{eqnarray}

\begin{Def}\hspace{-2mm}{\bf .}\label{LCmetrics}
Let  metrics $g$ and $\bg$ be 
 given by formulae (\ref{Canon_g}) and
(\ref{Canon_bg}) in a coordinate chart ${\cal U}$. Then we say that the metrics
$g$ and $\bg$ have {\em Levi-Civita local form
(of type m)}, and the coordinate    chart $\cal U$ is {\em a Levi-Civita}
coordinate chart (with respect to the  metrics). 
\end{Def}
Levi-Civita proved that the metrics $g$ and $\bg$  given by formulae
(\ref{Canon_g}) and (\ref{Canon_bg}) are geodesically equivalent.
If we replace 
$\phi_i$ by $\phi_i+c$, \  $ i=1,...,m$, where 
$c$ is a (positive for simplicity) constant, 
 in (\ref{Canon_g}) and (\ref{Canon_bg}), we  obtain
the following  one-parameter family of metrics,  
geodesically equivalent to $g$: 
\begin{equation} \nonumber
g_c(\dot{\bar x}, \dot{\bar x})=\frac{1}{(\phi_1+c)\cdots(\phi_m+c)}
\left\{\frac{1}{\phi_1+c}\Pi_1A_1+\cdots+\frac{1}{\phi_m+c}\Pi_mA_m\right\}.
\end{equation}

The next theorem is essentially due to Painlev\'e, see \cite{Levi-Civita}. 
\begin{Th}\hspace{-2mm}{\bf .}
If the metrics $g$ and $\bg$ are geodesically equivalent,
 then the function 
\begin{equation}\label{9} 
I_0\eqdef\left(\frac{{\rm det} (g)}{\det(\bg)}\right)^{\frac{2}{n+1}}\bg(\dot{\bar  x},\dot{\bar x}),
\end{equation}
is an integral of the geodesic flow of the 
metric $g$.
\end{Th}

\noindent Substituting $g_c$ instead of $\bar g$ in (\ref{9}),
 we 
obtain the following 
 one-parameter
family of  integrals
\begin{eqnarray}
I_c&\eqdef&\left(\frac{{\rm det}(g)}{{\rm det}(g_c)}\right)^{\frac{2}{n+1}}g_c(\dot{\bar x}, \dot{\bar x})=\nonumber\\
&=&C[(\phi_1+c)\cdots(\phi_m+c)]
\left\{\frac{1}{\phi_1+c}\Pi_1A_1+\cdots+\frac{1}{\phi_m+c}\Pi_mA_m\right\}\nonumber\\
&=&C\{L_1c^{m-1}+L_2c^{m-2}+\cdots+L_m\},\nonumber
\end{eqnarray}
\noindent where
\begin{eqnarray*}
L_1&=&\Pi_1A_1+\cdots+\Pi_mA_m,\;\mbox{ which is  twice the
 energy integral},\\
L_2&=&\si_1(\phi_2,...,\phi_m)\Pi_1A_1+\cdots+\si_1(\phi_1,...,\phi_{m-1})\Pi_mA_m,\\
L_3&=&\si_2(\phi_2,...,\phi_m)\Pi_1A_1+\cdots+\si_2(\phi_1,...,\phi_{m-1})\Pi_mA_m,\\
\vdots&&\\
L_m&=&(\phi_2...\phi_m)\Pi_1A_1+\cdots+(\phi_1...\phi_{m-1})\Pi_mA_m,
\end{eqnarray*}
 $\si_k$ denotes the elementary symmetric polynomial of degree $k$, 
and \\
$
C\eqdef\left[(\phi_1+c)^{k_1-1}\cdots(\phi_m+c)^{k_m-1}\right]^{\frac{2}{n+1}}
$
is a constant. Therefore 
the functions $L_k$, 
$k=1,...,m,$ are integrals of the geodesic flows of the metric $g$. 
We 
call these integrals  {\it Levi-Civita integrals}.

From the results of  \cite{Pain} it follows
 that Levi-Civita 
integrals are in involution.  More precisely, let $D=(d_j^i) $ 
be an  $m\times m$ matrix. 
Suppose that for any $i,j$ the element  $d_j^i$ depends only on  
the variables $\bar x_j$. Denote by $\Delta$ the determinant of the 
matrix $D$ and  by $\Delta_j^i$ the minor of the element $d_j^i$. 
In the paper \cite{Pain} it was  shown that, 
for arbitrary  functions 
$A_i(\bar x_i, \dot{\bar x}_i)$,  quadratic in velocities $\dot{\bar x}_i$, 
 the Lagrangian system with Lagrangian 
$$
T_1=\Delta\left( \frac{A_1(\bar x_1, \dot{\bar x}_1)}{\Delta_1^1 }+\frac{A_2(\bar x_2, \dot{\bar x}_2)}{\Delta^1_2 } + ... +\frac{A_m(\bar x_m, \dot{\bar x}_m)}{\Delta^1_m } \right) 
$$ admits $(m-1)$ integrals
$$
T_i=\Delta\left( A_1(\bar x_1, \dot{\bar x}_1)\frac{\Delta_1^i}{{(\Delta_1^1)^2} }+A_2(\bar x_2, \dot{\bar x}_2)\frac{\Delta^i_2}{(\Delta^1_2)^2 } + ... +A_m(\bar x_m, \dot{\bar x}_m)\frac{\Delta^i_m}{(\Delta^1_m)^2 
} \right),$$  where $i=2,...,m$, and if we identify 
 the tangent and cotangent bundles  the Lagrangian $T_1$ 
and consider the standard symplectic form on the cotangent bundle,
  then  the integrals are in involution.  

If we take $d_j^i=(\phi_j)^{m-i}$, then $\Delta$ and $\Delta_j^i$ are given by 
$$ \Delta^{i}_j= (-1)^{m-1}\sigma^{i-1}(\phi_1,\phi_2,...,\phi_{j-1}, \phi_{j+1},..., \phi_{m})\prod_{\alpha>\beta\ge 1,
 \alpha\ne j, \beta\ne j }(\phi_\alpha-\phi_\beta),$$
$$ \Delta= (-1)^{m}\prod_{\alpha>\beta\ge 1}(\phi_\alpha-\phi_\beta).$$     
Therefore, $$\frac{\Delta\Delta_j^i}{(\Delta_j^1)^2}= \sigma^{i-1}(\phi_1,\phi_2,...,\phi_{j-1}, \phi_{j+1},..., \phi_{m}) \Pi_j,$$ so  $T_i=-L_i$ and 
thus the integrals $L_i$ are in involution, q.~e.~d.

\section{Geodesic equivalence and  corresponding integrals}
\label{GeodEqInt} Let the
metrics $g$ and $\bar g$ on a
 manifold $M$ (of dimension n)  be  geodesically
equivalent.

\noindent Define
$$
U_g^rM\eqdef\left\{(x, \xi)\in{\cal T}M:\quad ||\xi||_g=r\right\},
$$
\noindent where 
$x\in M$,   $\xi\in {\cal T}_xM$
and $||\xi||_g\eqdef\sqrt{g(\xi,\xi)}=\sqrt{g_{ij}\xi^i\xi^j}$ is  the norm of the vector $\xi$ 
in the metric $g$.

By the geodesic flow of the metric $g$ we mean the Lagrangian 
system of differential equations 
$\frac{d}{dt}\left(\frac{\partial L}{\partial \dot x}\right)-\frac{\partial L}{\partial x}=0$ 
on ${\cal T}M$ with Lagrangian 
$L\eqdef\frac{1}{2}g_{ij}\dot x^i\dot x^j$.  
Because of the 
Legendre transformation,    the geodesic
 flow could be considered 
as a Hamiltonian system on ${\cal T}M$ (as a 
 symplectic 
form we take $\om_g\eqdef\dif[g_{ij}\xi^jdx^i]$) with the Hamiltonian 
$H_g\eqdef\frac{1}{2}g_{ij}\xi^i\xi^j$.

 Since the metrics $g, \bar g$ are geodesically  equivalent, 
the mapping
$\Phi:{\cal T}M\to {\cal T}M$,   
$\Phi(x,\xi)=\left(x, \xi\frac{||\xi||_g}{||\xi||_{\bar g}}\right)$,
takes the trajectories of the geodesic flow of the metric $g$ to the trajectories of the geodesic flow of the metric $\bar g$. 
This mapping  is a diffeomorphism 
(for $r\ne 0$),  takes $U_{g}^rM$ to $U_{\bar g}^rM$ and is 
trajectorial on  $U_{g}^rM$. Obviously the surfaces $U_g^r$, $U_{\bg}^r$ are
regular isoenergy surfaces $\{ H_g=\frac{r}{2}\}$,  
$\{  H_{\bar g}=\frac{r}{2}\}$.

By Theorem~\ref{MainTh},  in order to obtain   a 
family of
first integrals 
we have to find  the polynomial 
$\De^n(t)$ and divide
it by $(t-a)$. In our case 
 $H_g=H_{\bg}\circ\Phi$. Therefore  the  function $a$ from
Theorem \ref{MainTh}  equals to $\frac{||\xi||_{\bg}}{||\xi||_g}$.

\noindent In coordinates we have 
$$
\om_g=\dif[g_{ij}\xi^jdx^i]
$$
and
$$
\om_{\bg}=\dif[{\bg}_{ij}\xi^jdx^i].
$$
Therefore,
\begin{eqnarray*}
\Phi^{*}\om_{\bg}&=&\dif\left[\frac{||\xi||_g}{||\xi||_{\bg}}\bg_{ij}\xi^jdx^i\right]=\\
                 &=&\frac{\partial}{\partial x^k}
          \left[\frac{||\xi||_g}{||\xi||_{\bg}}\bg_{ij}\xi^j\right]dx^k\wedge dx^i-\\
                 &-&\frac{\partial}{\partial \xi^k}
          \left[\frac{||\xi||_g}{||\xi||_{\bg}}\bg_{ij}\xi^j\right]dx^i\wedge d\xi^k.\\
 \end{eqnarray*}
It is easy to see that at a point  $\xi\in{\cal T}_xM$ the quantities
$$
A_{ik}\eqdef-\frac{\partial}{\partial \xi^k}
          \left[\frac{||\xi||_g}{||\xi||_{\bg}}\bg_{ij}\xi^j\right]
$$
form an element of ${\cal T}_xM\otimes{\cal T}_xM$.
Without loss of generality we can assume that in  the space ${\cal T}_xM$ the metrics  $g$ and $\bg$ are given in principal axes. Then 
\begin{eqnarray*}
A_{ij}&\eqdef&-\rho^{i}(x)\frac{\partial}{\partial \xi^j}
\left(
\xi^i\frac{\sqrt{{\xi^1}^2+...+{\xi^n}^2}}{\sqrt{\rho^1{\xi^1}^2+...+\rho^n{\xi^n}^2}}
\right)=\\
             &=&\rho^i\de^i_j
\frac{||\xi||_g}{||\xi||_{\bg}}
              -\rho^i\xi^i\left(
   \frac{\frac{||\xi||_{\bg}}{||\xi||_g}-\rho^j\frac{||\xi||_g}{||\xi||_{\bg}}}
        {||\xi||_{\bg}^2}\xi^j\right)=\\
             &=&{\rm diag}(\mu_1,...,\mu_n)-A\otimes B.
\end{eqnarray*}
Here $\rho^i, i=1,...,n$ are common eigenvalues (here we allow $\rho^i$ to be equal to $\rho^j$  for some $i,j$) of the  metrics $g$ and $\bar g$, 
  $\mu_i\eqdef-\rho^i\frac{||\xi||_g}{||\xi||_{\bg}}$, $A_i\eqdef \rho^i\xi^i$
and $$B_i\eqdef\frac{\frac{||\xi||_{\bg}}{||\xi||_g}-\rho^i\frac{||\xi||_g}{||\xi||_{\bg}}}
        {||\xi||_{\bg}^2}\xi^i.$$

\noindent We have
$$
\begin{array}{rcl}
{\rm det}(\Phi^{*}\om_{\bg}-t\om_g)&={\rm det}&\left|
\begin{array}{c|c}
(*)&(A_{ij}+t\de_{ij})\\
\hline
-(A_{ij}+t\de_{ij})&0\
\end{array}
\right|\\
            &=&\det(A_{ij}+t\de_{ij})^2.
\end{array}
$$
\noindent Therefore,
\begin{equation}\nonumber
\Delta^n(t)=\det\left({\rm diag}(t+\mu_1,...,t+\mu_n)-a\otimes b\right).
\end{equation}

\begin{Lemma}\hspace{-2mm}{\bf .}\label{De} The following relation holds:
\begin{eqnarray}
\Delta^n(t)&=&(t+\mu_1)\cdots(t+\mu_n)-(a_1b_1)(t+\mu_2)\cdots(t+\mu_n)-\ldots\nonumber\\
           &-&(t+\mu_1)\cdots(t+\mu_{n-1})(a_nb_n).
\end{eqnarray}
\end{Lemma}

\noindent The lemma follows from induction considerations.

To divide the polynomial by  $(t-a)$ we shall use the 
Horner scheme. Suppose
that
$\De^n(t)=t^n+a_{n-1}t^{n-1}+\cdots+a_0$
and
$\de^{n-1}(t)=t^{n-1}+b_{n-2}t^{n-2}+\cdots+b_0$.
Then we have 
\begin{eqnarray}
&&b_{n-1}=a_n=1,\\
&&b_{n-2}=a_{n-1}+a\label{n-2},\\
&&\cdots\nonumber\\
&&b_{k}=a_{k+1}+ab_{k+1},\\
&&\cdots\nonumber\\
&&0=a_0+ab_0\label{-1}.
\end{eqnarray}
It follows from lemma \ref{De} that
\begin{eqnarray*}
a_0&=&(\mu_1...\mu_n)-(A_1B_1)(\mu_2...\mu_n)-\cdots
-(\mu_1...\mu_{n-1})A_nB_n=\\
   &=&(-1)^n\left(\frac{||\xi||_g}{||\xi||_{\bg}}\right)^n(\rho^1\cdots\rho^n).
\end{eqnarray*}
 Combining   with (\ref{-1}) we get 
$$
b_0=-\frac{a_0}{a}=
(-1)^{n+1}\left(\frac{||\xi||_g}{||\xi||_{\bg}}\right)^{n+1}(\rho^1\cdots\rho^n).
$$
\noindent Since $\frac{1}{2}g_{ij}\xi^i\xi^j$ is  an integral of the geodesic flow of the metric $g$, the function 
\begin{equation}\label{PainInt}
I_0\eqdef(\rho^1\cdots\rho^n)^{-\frac{2}{n+1}}\bg(\xi,\xi)
\end{equation}
\noindent is also an integral of the geodesic flow of the metric $g$. 
Using  Lemma \ref{De} we have 
\begin{eqnarray*}
a_{n-1}&=&(\mu_1+...+\mu_n)-(A_1B_1+...+A_nB_n)=\\
       &=&\frac{||\xi||_g}{||\xi||_{\bg}^3}\left\{({\rho^1}^2{\xi^1}^2+...+{\rho^n}^2{\xi^n}^2)-\right.\\
       &-&\left.(\rho^1+...+\rho^n)(\rho^1{\xi^1}^2+...+\rho^n{\xi^n}^2)\right\}-
\frac{||\xi||_{\bg}}{||\xi||_g}.
\end{eqnarray*}
\noindent Using (\ref{n-2}) we get 
\begin{eqnarray*}
b_{n-2}&=&a_{n-2}+a=\\
       &=&\frac{||\xi||_g}{||\xi||_{\bg}^3}\left\{({\rho^1}^2{\xi^1}^2+...+{\rho^n}^2{\xi^n}^2)-\right.\\
       &-&\left.(\rho^1+...+\rho^n)(\rho^1{\xi^1}^2+...+\rho^n{\xi^n}^2)\}\right.\\
\end{eqnarray*}
\noindent Therefore, the function
\begin{eqnarray*}
I_1&\eqdef&(\rho^1\cdots\rho^n)^{-\frac{3}{n+1}}
\left\{({\rho^1}^2{\xi^1}^2+...+{\rho^n}^2{\xi^n}^2)-\right.\\
&-&\left.(\rho^1+...+\rho^n)(\rho^1{\xi^1}^2+...+\rho^n{\xi^n}^2)\right\}
\end{eqnarray*}
\noindent is an integral. (It is easy to see  that
$
\frac{||\xi||^2_g}{||\xi||^2_{\bg}}=(\rho^1\cdots\rho^n)^{-\frac{2}{n+1}}
\frac{||\xi||_g^2}{I_0}.)$

Arguing as above, we see that the functions 
\begin{eqnarray*}
I_k&\eqdef&(\rho^1\cdots\rho^n)^{-\frac{k+2}{n+1}}
\left\{({\rho^1}^{k+1}{\xi^1}^2+...+{\rho^n}^{k+1}{\xi^n}^2)-\right.\\
&-&(\rho^1+...+\rho^n)({\rho^1}^k{\xi^1}^2+...+{\rho^n}^k{\xi^n}^2)+\cdots\\
&+&\left.(-1)^k\sigma_k(\rho^1,...,\rho^n)(\rho^1{\xi^1}^2+...+\rho^n{\xi^n}^2)\right\},
\end{eqnarray*}
are  integrals of the geodesic flow of the metric $g$, 
where by $\sigma_k$ we denote the elementary
symmetric polynomial of degree $k$.
It is obvious that $(-1)^k\sigma_k=c_k$ from Theorem \ref{Main_Topalov}, 
and therefore 
$I_k=\left(\frac{{\rm det}(g)}{{\rm det}(\bg)}\right)^{\frac{k+2}{n+1}}\bg(S_k\xi,\xi)$. Thus $I_k$, $k=0,..., n-1$,  are integrals of the geodesic flow of the metric
$g$, q.~e.~d.

\section{Liouville integrability\label{LiouvIntegr}}
\noindent 
The last step of  the proof of Theorem \ref{Main_Topalov}
is to 
verify  that the integrals  $I_0,...,I_{n-1}$  are in involution.
We proceed along  the following plan. 
First we show that it is sufficient to 
prove the involutivity in each Levi-Civita chart. Then we prove that in each Levi-Civita chart the integrals $I_0,..., I_{n-1}$ 
are linear combinations of Levi-Civita integrals, and therefore commute.

Let $g, \bg$ be metrics on $M$. 
A point $x\in M$ is called {\it stable}, 
if in 
 a neighborhood 
of $x$ the number 
of different   eigenvalues of the metrics 
$g, \bar g$ is   independent  of the  point.

Denote by ${\cal M}$ the set of stable 
 points of $M$.  The set ${\cal M}$ is an open subset of $M$. 
Obviously 
\begin{equation}\label{LCDecomp}
{\cal M}=\bigsqcup_{1\le q\le n}{\cal M}^q,
\end{equation}
where ${\cal M}^q$ denotes the set of stable
 points 
whose  number of distinct  common eigenvalues
equals  $q$. 
  Points $x\in M\setminus  {\cal M} $ 
are called  {\it  points of bifurcation}.

\begin{Lemma}\hspace{-2mm}{\bf .}\label{PolBifur}
The set ${\cal M}$ is everywhere dense in $M$. 
\end{Lemma}

{\bf Proof of Lemma \ref{PolBifur}.} 
Denote by $N(x)$ the number of  distinct  common eigenvalues 
of the metrics $g,\bar g$ at a  point $x$.
Recall that the  common eigenvalues 
of the metrics $g,\ \bar g$ at a
point $x\in M$ are roots of the
characteristic polynomial $P_x(t)={\rm  det}\left(G-tE\right)_{|x}$, where 
$G=\left(g^{i\alpha }\bar g_{\alpha j}\right)$. 
In particular, all roots of $P_x(t)$ are real.

Let us prove that, 
 for a sufficiently small neighborhood of an arbitrary point 
$x\in M$,   for any $y$ from the neighborhood 
the number $N(x)$ is no greater  than  $N(y)$.
Take a small $\epsilon>0$ and an arbitrary root $\rho$ of $P_x(t)$. 
Let us prove that for a sufficiently small neighborhood $U(x)\subset M$,  
for any $y\in U(x)$ 
there is a root $\rho_y$, $\rho-\epsilon <\rho_y<\rho+\epsilon$, 
 of the polynomial $P_y(t)$. 
If $\epsilon$ is  small,  then for a sufficiently small 
neighborhood  $U(x)$ of the point $x$, 
 for any $y\in U(x)$ 
the numbers $\rho+\epsilon$ and $\rho -\epsilon$ are not roots of $P_y(t)$. 
Consider the circle $S_\epsilon\eqdef \{z\in C:\quad  |z-\rho|=\epsilon\}$ on the complex plane $C$.  Clearly the number of roots (with multiplicities)  
of the polynomial $P_y$ inside the circle is equal to 
$$\frac{1}{2\pi i}\int_{S_\epsilon}\frac{P_y'(z)}{P_y(z)}dz.$$
Since for any $y\in U(x)$ there are no roots of $P_y$ on the circle $S_\epsilon$, 
then the function 
$$\frac{1}{2\pi i}\int_{S_\epsilon}\frac{P_y'(z)}{P_y(z)}dz
$$  continuously depends on  $y\in U(x)$, and therefore is a constant.
Clearly it is positive.  
Thus for any $y\in U(x) $ there is  at least one  root of $P_y$ that lies 
 between 
$\rho+\epsilon$ and $\rho -\epsilon$. Then for any $y$ from a sufficiently  
 small neighborhood of $x$ we have 
$N(y)\ge N(x)$.

Now let us prove the lemma. Evidently the set ${\cal M}$
 is an open subset of $M$. Then it is sufficient to prove that for any 
open  subset $U\subset  M$  there is a stable point $x\in U$. Suppose otherwise,
 i.e.
let all the points of  $U$ be  points of bifurcation.
Take a point $y\in M$ with maximal value of the function $N$ on   it. 
 We have that in a  neighborhood $U(y)$  of the point $y$ the function 
$N$ is constant and equals  $N(y)$. Then  the point $y$ is a stable point, 
and  we get  
a contradiction,  q.~e.~d.

Now let the metrics $g , \bar g$ be geodesically equivalent. 
Since the set of points of 
bifurcation is nowhere dense, it is sufficient to prove the involutivity 
in each Levi-Civita  chart.  
Let the metrics $g$ and $\bg$ be given by 
\begin{eqnarray}
g(\dot{\bar x}, \dot{\bar x})  &=&\Pi_1(\bar x)A_1(\bar x_1,\dot{\bar x}_1)+
\Pi_2(\bar x)A_2(\bar x_2,\dot{\bar x}_2)+\cdots+\nonumber\\
&+&\Pi_m(\bar x)A_m(\bar x_m,\dot{\bar x}_m),\\
\bg(\dot{\bar x}, \dot{\bar x}) &=&\rho^1\Pi_1(\bar x)A_1(\bar x_1,\dot{\bar x}_1)+
\rho^2\Pi_2(\bar x)A_2(\bar x_2,\dot{\bar x}_2)+\cdots\nonumber\\
&+&\rho^m\Pi_m(\bar x)A_m(\bar x_m,\dot{\bar x}_m).
\end{eqnarray}
\noindent We show that the integrals $I_k$ are linear combinations of the 
Levi-Civita integrals. 
  We have 
\begin{equation}
\bar G={\rm diag}(\underbrace{\rho^1,...,\rho^1}_{k_1},...,
\underbrace{\rho^m,...,\rho^m}_{k_m}),
\end{equation}
\noindent where $\rho^k=\frac{1}{(\phi_1...\phi_m)}\frac{1}{\phi_k}$.
It is easy to check that
\begin{equation}
S_k=(-1)^k{\rm diag}(\underbrace{\si_k^1,...,\si_k^1}_{k_1},...,
\underbrace{\si_k^m,...,\si_k^m}_{k_m}),
\end{equation}
\noindent where
\begin{equation}
\si_k^l\eqdef\sigma_k(\underbrace{\rho^1,...,\rho^1}_{k_1},...,
\underbrace{\rho^l,...,\rho^l}_{k_l-1},...,
\underbrace{\rho^m,...,\rho^m}_{k_m}).
\end{equation}
We have
\begin{eqnarray}
\si_k^1&=&\frac{1}{(\phi_1...\phi_m)^k}
\sigma_k\Bigl(\underbrace{\frac{1}{\phi_1},...,\frac{1}{\phi_1}}_{k_1-1},...,
\underbrace{\frac{1}{\phi_m},...,\frac{1}{\phi_m}}_{k_m}\Bigr)=\\
&=&\frac{1}{(\phi_1...\phi_m)^k}\sum_{|\alpha|=k}
{k_1-1\choose\alpha_1}{k_2\choose\alpha_2}\cdots{k_m\choose\alpha_m}
\frac{1}{\phi_1^{\alpha_1}}\frac{1}{\phi_2^{\alpha_2}}\cdots\frac{1}{\phi_m^{\alpha_m}},\\
\end{eqnarray}
\noindent where $|\alpha|\eqdef\alpha_1+\cdots+\alpha_m$ and $\alpha_i\ge 0$.
Substituting   ${k_l-1\choose\alpha_l}+{k_l-1\choose\alpha_l-1}$ for
  ${k_l\choose\alpha_l}$  
 (we assume that
${k\choose 0}=1$, ${k\choose -1}=0$, $k\ge 0$)
for $2\le l\le m$ we obtain
\begin{eqnarray*}
\si_k^1&=&\frac{1}{(\phi_1...\phi_m)^k}\left(B_k+
B_{k-1}\si_1\Bigl(\frac{1}{\phi_2},...,\frac{1}{\phi_m}\Bigr)+\cdots+\right.\\
&+&\left.B_{k-m+1}\sigma_{m-1}\Bigl(\frac{1}{\phi_2},...,\frac{1}{\phi_m}\Bigr)
\right),
\end{eqnarray*}
\noindent where
\begin{equation}
B_k\eqdef\sum_{|\alpha|=k}{k_1-1\choose\alpha_1}\cdots{k_m-1\choose\alpha_m}
\frac{1}{\phi_1^{\alpha_1}}\cdots\frac{1}{\phi_m^{\alpha_m}}.\\
\end{equation}
\noindent Note that
\begin{equation}
\left(\frac{{\rm det} (g)}{{\rm det} (\bar g)}\right)^{\frac{k+2}{n+1}}=C_k(\phi_1...\phi_m)^{k+2},
\end{equation}
\noindent where
$C_k=[\phi_1^{k_1-1}...\phi_m^{k_m-1}]^{\frac{k+2}{n+1}}$.
\noindent Therefore,
\begin{eqnarray}
I_k&\eqdef&\left(\frac{{\rm det}(g)}{{\rm det}(\bg)}\right)^{\frac{k+2}{n+1}}\bg(S_k\dot{\bar x},\dot{\bar x})=\nonumber\\
&=&(-1)^kC_k(\phi_1...\phi_m)^{k+2}\left\{\rho^1\si_k^1\Pi_1A_1+\cdots+\rho^m\si_k^m
\Pi_m A_m\right\}=\nonumber\\
&=&(-1)^kC_k(\phi_1...\phi_m)^{k+2}
\left\{\frac{1}{\phi_1...\phi_m}\frac{1}{\phi_1}\left\{
\frac{1}{(\phi_1...\phi_m)^k}\left(B_k+\right.\right.\right.\nonumber\\
&+&\left.\left.\left.\cdots+B_{k-m+1}\sigma_{m-1}\Bigl(\frac{1}{\phi_2},...,\frac{1}{\phi_m}\Bigr)
\right)\right\}\Pi_1A_1+\cdots\right\}=\nonumber\\
&=&(-1)^kC_k\left\{B_kL_m+B_{k-1}L_{m-1}+\cdots+B_{k-m+1}L_1\right\},\label{Decomp}
\end{eqnarray}
\noindent where $L_i$ are Levi-Civita integrals.

Finally, since the integrals $I_0,..., I_{n-1}$ are linear combinations of  
Levi-Civita integrals with constant coefficients,
 and since Levi-Civita  
integrals commute, then the integrals $I_0,..., I_{n-1}$ also commute, q.~e.~d.

\begin{Rem}\hspace{-2mm}{\bf .}\label{rank}
Let $m$ be the number of distinct common 
eigenvalues of  geodesically equivalent metrics $g, \bar g$ at
  a point $x$. Then in a neighborhood $U$ of  the point 
 $x$  the number of functionally independent almost everywhere 
  Levi-Civita integrals is  no less than  $m$. Therefore the dimension of the space 
generated by the differentials 
$(dI_0, dI_1,..., dI_{n-1})$ no less than  $m$ at
almost all points    of ${\cal T}U$. 
\end{Rem}

\section{Topological obstructions}\label{Obstructions}
 Corollary \ref{T} follows immediately  from  the following theorem. 
Recall that a  group  $G$ is {\it almost commutative}, if there exists 
 a commutative subgroup $P\subset G$ of finite index. 

\begin{Th}[Taimanov, \cite{Taiman}]\hspace{-2mm}{\bf .} \label{Taim} 
If  a real-analytic closed manifold $M^n$ with a real-analytic
metric satisfies at least one of the conditions:
\begin{itemize}
\item[a)] $\pi_1(M^n)$ is not almost commutative
\item[b)] $dim H_1(M^n;{\bf Q})>dim M^n$, 
\end{itemize}
then the geodesic flow on $M^n$ is not analytically integrable. 
\end{Th}

\noindent{\bf Proof of Corollary \ref{T}.}
If metrics $g, \bar g$ are real-analytic  and geodesically equivalent,
then the integrals $I_0,...,I_{n-1}$ are also real-analytic. 
If the 
metrics are strictly non-proportional at least at   one point 
of $M^n$, then 
the integrals  are functionally independent almost everywhere in a neighborhood of that
 point. Since the integrals are real-analytic,  then 
they are functionally independent almost everywhere  and  
  we can apply Theorem~\ref{Taim}, q.~e.~d. 

\vspace{.2cm}

\noindent{\bf Proof of Corollary \ref{dim2}.} Let metrics $g,\bar g$ on a surface $M^2$ be geodesically equivalent. Using Theorem~\ref{Main_Topalov} we have that the function 
 $I_0=\left(\frac{{\rm det}(g)}{{\rm det}(\bar g)}\right)^{\frac{2}{n+1}}\bar g(\xi, \xi)$ is an integral of the geodesic flow of the metric $g$. In one direction Corollary~\ref{dim2}
is proved. In other direction the statement of Corollary~\ref{dim2}
can be verified by direct calculation,
 and it was done in \cite{Whi}. 

\vspace{2mm}
 \noindent{\bf Proof of Corollaries \ref{Kol}, \ref{torus}, \ref{sphere}, \ref{nontrivial}.}  Let $g$ be 
a metric on a surface $M^2$. 
The following lemma is essentially due to \cite{Birkhoff}, see also \cite{Kol}.
For simplicity assume that the surface $M^2$  is oriented, otherwise 
finitely cover the surface by an 
 oriented one. 
Consider the complex structure on $M^2$ corresponding to the metric $g$. 
Let $z$ be a complex coordinate in a open domain $U\subset M^2$.
Consider the complex momentum $p$.  We shall denote by $\bar z$ and  
$\bar p$ the 
complex conjugation of $z$ and  $p$ respectively. 
In complex variables the Hamiltonian $H:{\cal T}^*M^2\to R$
of the geodesic flow of the metric $g$ reads $\frac{p\bar p}{\lambda(z)}$, 
where $\lambda $ is a real-valued function.
 Suppose that the  
real-valued  function 
$$ 
  F=A(z)p^2+B(z)p\bar p+\bar A(z)\bar p^2
$$
is an integral of the geodesic flow of the metric $g$.

\begin{Lemma}\hspace{-2mm}{\bf .}\label{Birkhoff}
The form $\frac{1}{A(z)}dzdz$ is meromorphic. 
\end{Lemma}

\begin{Rem}\hspace{-2mm}{\bf .}
If  the Hamiltonian and  the integral are proportional at  each point of $M^2$, i.e. if $F\equiv \alpha(z)H$, where $\alpha:M^2\to R$,
  then by definition 
put  
$\frac{1}{A(z)}dzdz$ equal  zero. 
\end{Rem}

\noindent {\bf Proof of Lemma~\ref{Birkhoff}}. 
Since $F$ is an integral of the Hamiltonian system with the Hamiltonian
$H$, the  Poisson bracket $\{H, F\}$ equals zero. We have 
\begin{equation}\label{zero}
\{H, F\}= H_pF_z-H_zF_p+ H_{\bar p}F_{\bar z}-H_{\bar z}F_{\bar p}=0
\end{equation}
\noindent On the right side  of (\ref{zero}) each term is 
a polynomial of third degree in  momenta. Then the bracket is also a 
polynomial of third degree in momenta. In order for  a
polynomial  to equal zero, 
all coefficients must be zero, 
in particular the coefficient of  $p^3$. 
Thus $\frac{A_{\bar z}}{\lambda}$ equals zero, 
and $A$ is holomorphic. Then $\frac{1}{A(z)}$ is meromorphic, q.~e.~d.

Let $g, \bar g$ be geodesically equivalent
 metrics 
on a closed surface $M^2$ of Euler characteristic $\chi(M^2)$. Then the function 
$I_0=\left(\frac{{\rm det}(g)}{{\rm det}(\bar g)}\right)^{\frac{2}{n+1}}\bar g(\xi, \xi)$ 
is an integral of the geodesic flow of the metric $g$, 
and is quadratic in momenta (if we identify with the help of the metric $g$
 the tangent and  cotangent bundles  of  $M^2$).
 Consider the form  
$\frac{1}{A(z)}dzdz$ corresponding to the integral $I_0$. 
Suppose that the form is not identical zero.  
 For a meromorphic  2-form on a closed Riemann surface, the  number 
of poles $P$   minus the number of zeros $Z$ 
is equal to  twice the   Euler characteristic. 
It is easy to 
see that  the form $\frac{1}{A(z)}dzdz$
has 
no zeros  (otherwise the metric $\bar g$ has singularities).
Then  $P=2\chi(M^2)$, 
and the Euler characteristic $\chi(M^2)$ can not 
be negative, q.~e.~d. 
Now assume the metrics are proportional at each point of an open subset $U\subset M^2$. Since the form is meromorphic, it must be zero.  Thus 
$\bar g=\alpha(z)g$, where 
$\alpha$ is a function on $M^2$. Let us   show 
 that  the function $\alpha$ is constant. Actually, 
 $I_0=2\left(\frac{1}{\alpha}\right)^{\frac{1}{3}}H$ 
(here we identify ${\cal T}^*M$ and ${\cal T}M$ with the help of the metric $g$). We have 
$$\{H,I_0\}=\{H, 2\left(\frac{1}{\alpha}\right)^{\frac{1}{3}}H\}=
\{H, H\}2\left(\frac{1}{\alpha}\right)^{\frac{1}{3}}+ 2H\{\left(\frac{1}{\alpha}\right)^{\frac{1}{3}},H\}. 
$$ 
Since $\{ H, H\}$ equals zero,
we have that
$\{\left(\frac{1}{\alpha}\right)^{\frac{1}{3}},H\}$
 equals zero and the function $\alpha$ is constant.
 This proves Corollaries~\ref{Kol},\ref{nontrivial}.

\begin{Rem}\hspace{-2mm}{\bf .}
For non-orientable surfaces the sign of the Euler characteristic
coincides with the sign of the Euler characteristic
  of the oriented covering. 
Therefore Corollary~\ref{Kol} is true also for 
   non-orientable surfaces.
\end{Rem}
   
It is easy to 
see that  the form $\frac{1}{A(z)}dzdz$
 has poles precisely  at 
 points,  where the metrics are proportional. 
If the surface $M^2$ is the torus, 
then $\chi(M^2)=0$ 
and  either the metrics $g, \bar g$ are  proportional 
at  every point, or there are no points of proportionality of the
 metrics. This proves Corollary~\ref{torus}.

The following lemma is essentially due to Kolokol'tzov  
\cite{Kol}. It completes the proof of Corollary~\ref{sphere}. 

\begin{Lemma}\hspace{-2mm}{\bf .}
On the sphere $S^2$ there are  the following 
three possibilities for the form $\frac{1}{A(z)}dzdz$.  
\begin{itemize}
\item[1.] The form $\frac{1}{A(z)}dzdz $ is identical zero. 
\item[2.] The form  $\frac{1}{A(z)}dzdz $ has exactly two zeros (both zeros are of multiplicity two).
\item[3.] The form  $\frac{1}{A(z)}dzdz $ has exactly four zeros.
 \end{itemize}
In the second case the  metric $g$ admits a non-trivial 
Killing vector field. 
\end{Lemma}

\noindent {\bf Proof of Corollary~\ref{killing}.} 
Because of  Noether's  theorem, 
if a metric admits a (non-trivial) Killing
vector field, 
then the geodesic flow of the metric admits a (non-trivial) 
 integral, linear in velocities,  and vice versa. 
     
Suppose the 
 function 
$$F_1=\sum^n_{i=1}a_i(x)\xi^i$$
is constant on 
the trajectories of the geodesic flow of the metric
$\bar g$. Then  the function
$$\Phi^*F_1=
\frac{||\xi||_{g}}{||\xi||_{\bar g}}\sum^n_{i=1}a_i(x)\xi^i$$ is 
constant on the trajectories of the geodesic flow of the metric $g$. Since the function 
$
I_0=\left(\frac{{\rm det}(g)}{{\rm det}(\bar g)}\right)^{\frac{2}{n+1}}\bar g(\xi, \xi)
$
is an integral  of the geodesic flow of the metric $g$, and since the function 
$||\xi||_{g}=\sqrt{g(\xi,\xi)}$ is also an integral 
of the geodesic flow of the metric $g$, then the 
function  
$$\frac{\sqrt{g(\xi,\xi)}}{\sqrt{I_0}}\Phi^*F_1=
\left(\frac{{\rm det}(g)}{{\rm det}(\bar g)}\right)^{\frac{1}{n+1}}\sum^n_{i=1}a_i(x)\xi^i, 
$$ linear in velocities,  
is also an integral 
of the geodesic flow of the metric $g$, q.~e.~d.

\section{Geodesically equivalent metrics on the ellipsoid.   }\label{Ellipsoid}
{\bf Proof of Theorem~\ref{ellipsoid}. }
We show that in the elliptic  coordinate system the restriction of 
the metrics 
$$
ds^2\eqdef\sum_{i=1}^n(dx^i)^2 \ {\rm  and} \  \ dr^2\eqdef\frac{1}{\sum_{i=1}^{n} \left(\frac{x^i}{a_i}\right)^2}\left(\sum_{i=1}^n\frac{(dx^i)^2}{a_i}\right) 
$$
to  
the ellipsoid 
  $\sum_{i=1}^n\frac{(x^i)^2}{a_i}=1$ have Levi-Civita local form, and therefore are geodesically equivalent. 
More precisely, 
consider elliptic coordinates $\nu^1, ..., \nu^n$. 
Without loss of generality we can assume that $a^1<a^2<...<a^n$. Then 
the relation between 
the elliptic 
 coordinates  $\bar\nu$ and  the Cartesian coordinates $\bar x$ is 
given by 
\begin{equation}\nonumber 
x^i=\sqrt{\frac{\prod_{j=1}^n(a^i-\nu^j)}{\prod_{j=1, j\ne i}^n(a^i-a^j)}} . 
\end{equation}
Recall that the elliptic coordinates are non-degenerate almost everywhere, 
and the set $$\{\nu^1=0, a_1<\nu^2<a_2, a_2<\nu^3<a_3,..., a_{n-1}<\nu^{n}<a^n\}$$     
is the part of the ellipsoid,  lying in the
 quadrant $\{x^1>0, x^2>0, ... , x^n>0\}$.
Since for any $i$ the symmetry $x^i\to -x^i$
takes the ellipsoid to the ellipsoid and preserves  the metrics $ds^2$ and  $dr^2$, it is sufficent to check the statement of the theorem only in  the   quadrant $\{x^1>0, x^2>0, ... , x^n>0\}$. 

In the elliptic coordinates  
 the restriction of the metric $ds^2$ to the ellipsoid 
has the following form 
\begin{equation}\nonumber
\sum_{i=1}^n\Pi_iA_i(d\nu^i)^2, 
\end{equation}
where $\Pi_i\eqdef\prod_{j=1, j\ne i}^n(\nu^i-\nu^j)$,
 and $A_i\eqdef\frac{\nu^i}{\prod_{j=1}^n(a^j-\nu^i) }$. 
The restriction of the metric $dr^2$ to the ellipsoid 
is  
\begin{equation}\nonumber
(a^1a^2...a^n)\sum_{i=1}^n\rho^i\Pi_iA_i(d\nu^i)^2, 
\end{equation}
where $\rho^i\eqdef \frac{1}{\nu^i(\nu^1\nu^2...\nu^n)}.$
We see that the metrics $ds^2$, 
$dr^2$ have Levi-Civita local form, and therefore are geodesically 
equivalent, q.~e.~d.

\end{document}